\documentclass{amsart}
\usepackage{amssymb,latexsym,amsmath,amscd,graphicx,graphics,epic,eepic,enumerate}

\theoremstyle{plain}
\newtheorem{theorem}{Theorem}[section]

\newtheorem{proposition}[theorem]{Proposition}

\theoremstyle{definition}

\theoremstyle{remark}
\newtheorem{remark}[theorem]{Remark}

\newcommand{\nat}{\mathbb{N}}
\newcommand{\zed}{\mathbb{Z}}

\newcommand{\C}{\mathbb{C}}
\newcommand{\fil}{\mathcal{F}}

\newcommand{\gpm}{g_p^{\min}}
\newcommand{\gpM}{g_p^{\max}}
\newcommand{\gnm}{g_n^{\min}}
\newcommand{\gnM}{g_n^{\max}}
\newcommand{\la}{\mathcal{L}}

\begin{document}

\title{The Khovanov-Rozansky Cohomology and Bennequin Inequalities}

\author{Hao Wu}

\address{Department of mathematics and Statistics, Lederle Graduate Research Tower, 710 North Pleasant Street, University of Massachusetts, Amherst, MA 01003-9305, USA}

\email{wu@math.umass.edu}



\begin{abstract}
We review Bennequin type inequalities established using various versions of the Khovanov-Rozansky cohomology. Then we give a new proof of a Bennequin type inequality established by the author \cite{Wu5}, and derive new Bennequin type inequalities for knots using Gornik's version of the Khovanov-Rozansky cohomology, which generalize those established by Shumakovitch \cite{Shu}, Plamenevskaya \cite{Pl4} and Kawamura \cite{Kawa} using the Rasmussen invariant.
\end{abstract}

\maketitle

\section{Introduction}

In 1983, Bennequin \cite{Ben} proved that
\begin{equation}\label{Bi}
\chi(B) \leq w+O
\end{equation}
for any closed braid $B$, where $\chi(B)$ is the maximal Euler characteristic of a Seifert surface of $B$,  $O$ is the number of Seifert circles of $B$, and $w$ is the writhe of $B$. This inequality is now called the Bennequin inequality. Rudolph \cite{Ru} refined the Bennequin inequality and proved the slice Bennequin inequality
\begin{equation}\label{sBi}
\chi_s(B) \leq w+O,
\end{equation}
where $\chi_s(B)$ is defined to be the maximal Euler characteristic, $\chi(F)$, of an oriented smoothly embedded compact surface $F$ without closed components in $D^4$ bounded by $B$. Using the HOMFLY polynomial, Franks and Williams \cite{FW} and Morton \cite{Mo} established the following Morton-Franks-Williams inequality,
\begin{equation}\label{HOMFLYBi}
w-O ~\leq~ \min-\deg_a P_B ~\leq~ \max-\deg_a P_B ~\leq~ w+O,
\end{equation}
where $P_B$ is the HOMFLY polynomial of $B$ and $a$ is the framing variable of $P_B$.

Although these inequalities were originally proved only for closed braids, from the Corollary in Section 4 of \cite{Ya}, it is easy to see that these are actually true for all link diagrams. In this paper, we call these and other similar inequalities the Bennequin type inequalities. Ferrand \cite{Fer} gave a detailed account of developments of Bennequin type inequalities up to year 2000.

In 1999, Khovanov \cite{K1} defined a new link invariant categorifying the Jones polynomial, which is now called the Khovanov homology. Lee \cite{Lee2} constructed a perturbed version of the Khovanov homology, and showed its invariance. Rasmussen \cite{Ras1} used Lee's construction to defined the Rasmussen invariants $s$ for knots, and showed that the absolute value of $s$ is bounded from above by twice the slice genus of the knot. In 2004, Shumakovitch \cite{Shu} and Plamenevskaya \cite{Pl4} proved that, for any diagram of a knot $K$,
\begin{equation}\label{SPBi}
w-O\leq s(K)-1,
\end{equation}
which is sharper than inequality \eqref{sBi} for knots. Kawamura \cite{Kawa} refined this result, and proved that, for any diagram of a knot $K$,
\begin{equation}\label{KBi}
w-O_\geq+O_< ~\leq~ s(K)-1,
\end{equation}
where $O\geq$ and $O_>$ are defined below.

In 2004, Khovanov and Rozansky \cite{KR1} generalized the Khovanov homology and defined a link invariant categorifying the $SO(n)$-HOMFLY polynomial, which is now called the Khovanov-Rozansky $sl(n)$-cohomology. Using this invariant, the author \cite{Wu5} proved that, for any diagram of a link $L$,
\begin{equation}\label{aKRBi}
w-O ~\leq~ \liminf_{n\rightarrow+\infty} \frac{\gnm(L)}{n-1} ~\leq~ \limsup_{n\rightarrow+\infty} \frac{\gnM(L)}{n-1} ~\leq~ w+O,
\end{equation}
where $\gnm$ and $\gnM$ are defined below. This inequality is sharper than inequalities \eqref{sBi} and \eqref{HOMFLYBi}. See \cite{Wu5,Wu7} for more details.

Shortly after Khovanov and Rozansky \cite{KR1} defined their $sl(n)$-cohomology, Gornik \cite{Gornik} generalized Lee's construction and gave a perturbed version of the Khovanov-Rozansky $sl(n)$-cohomology, which was recently proved to be a link invariant by the author \cite{Wu7}.

The goal of the present paper is to give a new proof of inequality \eqref{aKRBi} and establish new Bennequin type inequalities using Gornik's version of the Khovanov-Rozansky $sl(n)$-cohomology, which generalize the inequalities \eqref{SPBi}, \eqref{KBi}.

In the rest of this paper, we fix an integer $n\geq2$ and $p(x)=x^{n+1}-(n+1)x$, and use the terminologies introduced in \cite{Wu7}. In particular, write $H_n=H_{x^{n+1}}$ for the original Khovanov-Rozansky $sl(n)$-cohomology defined in \cite{KR1}, and $H_p$ for the Khovanov-Rozansky cohomology defined using $p(x)$, which is the perturbed version defined by Gornik \cite{Gornik}. We also define
\begin{eqnarray*}
\gpM & = & \max\{k~|~\fil^kH_p/\fil^{k-1}H_p\neq0\}, \\
\gpm & = & \min\{k~|~\fil^kH_p/\fil^{k-1}H_p\neq0\}, \\
\gnM & = & \max\{k~|~H_n^k\neq0\}, \\
\gnm & = & \min\{k~|~H_n^k\neq0\},
\end{eqnarray*}
where $\fil$ is the quantum filtration of $H_p$, and $H_n^k$ is the subspace of $H_n$ consisting of elements with quantum grading $k$. These are numerical invariants for links.

The next two theorems serve as the starting point of our argument.

\begin{theorem}[\cite{Wu7}]\label{slice-genus}
For any knot $K$ in $S^3$,
\begin{eqnarray*}
(n-1)\chi_s(K) ~ \leq & \gpM(K) & \leq~ (n-1)(2-\chi_s(K)), \\
(n-1)(\chi_s(K)-2) ~\leq & \gpm(K) & \leq~ -(n-1)\chi_s(K),
\end{eqnarray*}
where $\chi_s(K)$ is the slice Euler characteristic of $K$, which is defined to be the maximal Euler characteristic, $\chi(S)$, of an oriented smoothly embedded compact surface $S$ without closed components in $D^4$ bounded by $K$.
\end{theorem}

\begin{theorem}[\cite{Wu5}]\label{res-braid}
Let $B$ be a closed braid with $O$ strands, $c_+$ positive crossings and $c_-$ negative crossings. Then
\[
(n-1)(w-O)-2c_- ~\leq~ \gnm(B) ~\leq~ \gnM(B) ~\leq~ (n-1)(w+O)+2c_+,
\]
where $w=c_+-c_-$ is the writhe of $B$.
\end{theorem}

Theorem \ref{slice-genus} is proved at the end of \cite{Wu7} using link cobodisms. Theorem \ref{res-braid} is proved in \cite{Wu5} by inducting on the "weight" of a resolved closed braid. Clearly, inequality \eqref{aKRBi} is an immediate consequence of Theorem \ref{res-braid}. In Section \ref{res-braid-sec}, we give a new proof of Theorem \ref{res-braid} using the composition product recently established by Wagner \cite{Wag}.

Combining Theorems \ref{slice-genus} and \ref{res-braid}, we determine $\gpM$ for negative closed braids and $\gpm$ for positive closed braids. Then we do an induction using some simple link cobodisms, and prove the following Bennequin type inequality for knots.

\begin{theorem}\label{bennequin}
For any knot $K$ in $S^3$ and any knot diagram $D_K$ of $K$,
\[
(n-1)(w-O) ~\leq~ \gpm(K) ~\leq~ \gpM(K) ~\leq~ (n-1)(w+O),
\]
where $w$ is the writhe of $D_K$ and $O$ is the number of Seifert circles of $D_K$.
\end{theorem}

Note that, when $n=2$, the Rasmussen invariant of a knot $K$ is $s(K)= \gpM(K)-1=\gpm(K)+1$. So Theorem \ref{bennequin} generalizes inequality \eqref{SPBi}.

In the setup of Theorem \ref{bennequin}, let $O_>$ (resp. $O_<$) be the number of Seifert circles of $D_K$ adjacent to only positive (resp. negative) crossings, and $O_\leq = O-O_>$ (resp. $O_\geq =O-O_<$). Applying an argument by Rudolph \cite{Ru1}, we refine Theorem \ref{bennequin}, and prove the following sharper Bennequin inequality, which generalizes inequality \eqref{KBi}.

\begin{theorem}\label{sharper-bennequin}
Let $K$ be any knot, and $D_K$ any diagram of $K$.
\begin{enumerate}[(a)]
    \item If $D_K$ has only positive crossings, then
          \[
          (n-1)(w-O) ~=~ \gpm(K) ~\leq~ \gpM(K) ~\leq~ (n-1)(w-O+2).
          \]
    \item If $D_K$ has only negative crossings, then
          \[
          (n-1)(w+O-2) ~\leq~ \gpm(K) ~\leq~ \gpM(K) ~=~ (n-1)(w+O).
          \]
    \item If $D_K$ has both positive and negative crossings, then
          \[
          (n-1)(w-O_\geq+O_<) ~\leq~ \gpm(K) ~\leq~ \gpM(K) ~\leq~ (n-1)(w+O_\leq-O_>).
          \]
\end{enumerate}
\end{theorem}

\begin{remark}
(a) Note that all the terms in Theorems \ref{bennequin} and \ref{sharper-bennequin} are defined for links. So, although the proofs of these two theorems in the present paper work only for knots, it is likely that one can generalize these inequalities to links, at least in some weaker form.

(b) Other interesting properties of Lee's homology and, in particular, the Rasmussen invariant have been found since the definition of these invariants. Now the invariance of Gornik's version of the Khovanov-Rozansky cohomology is established, and the Rasmussen invariant is generalized to $\forall~ n\geq 2$ (See \cite{Wu7}). It is natural to ask if these properties of Lee's homology and the Rasmussen invariant can be generalized. For example, is it true that $\gpM-\gpm=2(n-1)$ or, at least, $\gpM-\gpm\leq2(n-1)$  for all knots? Is the generalized Rasmussen invariant a concordance invariant for knots?

(c) The Bennequin type inequalities above induce upper bounds for the self-linking number of a transversal link (or knot) in the standard contact $S^3$, which, in turn, induce upper bounds for the Thurston-Bennequin number of a Legendrian link (or knot) in the standard contact $S^3$. There are upper bounds for the Thurston-Bennequin number that are not obtained this way. The Kauffman polynomial provides one such upper bound. (See, e.g. \cite{Fer,FT}.) Ng \cite{NgKB} also established such an upper bound using the Khovanov homology. It seems unlikely that one can generalize Ng's result to the Khovanov-Rozansky cohomology. Maybe the generalization will only come after the categorification of the Kauffman polynomial is successful.

(d) It is interesting to compare the Bennequin type inequalities mentioned in the present paper to the one established by Plamenevskaya \cite{Pl2} using the knot concordance invariant $\tau$ from the Heegaard-Floer homology \cite{OS4}.
\end{remark}

\section{Composition product and the proof of Theorem \ref{res-braid}}\label{res-braid-sec}

A planar diagram $\Gamma$ is called a graph if it has the following
properties:
\begin{enumerate}[(a)]
    \item $\Gamma$ consists of two types of edges: regular edges
    and wide edges. These edges intersect only at their endpoints.
    \item Regular edges are disjoint from each other.
    Wide edges are disjoint from each other.
    \item Each regular edge is oriented, and contains at least one
    marked point. Open endpoints of regular edges are marked.
    \item Each wide edge has exactly two regular edges entering at one endpoint,
    and exactly two regular edges exiting from the other endpoint.
\end{enumerate}
In the rest of this section, all graphs are closed, i.e., have no open endpoints.

Let us recall the composition product established by Wagner
\cite{Wag}. A labeling of a closed graph $\Gamma$ is a function from
the set of all regular edges of $\Gamma$ to $\{1,2\}$, which
satisfies that, at each wide edge $E$, the number of adjacent
regular edges labeled by $1$ (resp. $2$) entering $E$ is equal to
the number of adjacent regular edges labeled by $1$ (resp. $2$)
exiting $E$. Let $\la(\Gamma)$ be the set of all labelings of
$\Gamma$.

\begin{figure}[ht]

\setlength{\unitlength}{1pt}

\begin{picture}(360,60)(-180,40)


\put(-52.5,85){\vector(-1,1){15}}

\put(-47.5,85){\vector(1,1){15}}

\put(-52.5,75){\line(-1,-1){15}}

\put(-47.5,75){\line(1,-1){15}}

\put(-75,95){\small{$2$}}

\put(-75,60){\small{$2$}}

\put(-30,95){\small{$1$}}

\put(-30,60){\small{$1$}}

\put(-75,40){$\left\langle E|\Gamma|f \right\rangle=1$}


\put(52.5,85){\vector(1,1){15}}

\put(47.5,85){\vector(-1,1){15}}

\put(52.5,75){\line(1,-1){15}}

\put(47.5,75){\line(-1,-1){15}}

\put(70,95){\small{$2$}}

\put(70,60){\small{$2$}}

\put(25,95){\small{$1$}}

\put(25,60){\small{$1$}}

\put(25,40){$\left\langle E|\Gamma|f \right\rangle=-1$}

\linethickness{5pt}

\put(-50,75){\line(0,1){10}}

\put(50,75){\line(0,1){10}}

\end{picture}

\caption{Non-zero local interactions}\label{interaction}

\end{figure}
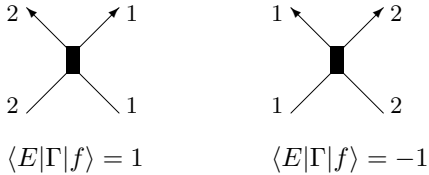

For a graph $\Gamma$, a labeling $f\in\la(\Gamma)$ and a wide edge $E$ of $\Gamma$, define the local interaction $\left\langle E|\Gamma|f \right\rangle$ to be $0$ except in the two cases depicted in Figure \ref{interaction}. And define
\[
\left\langle \Gamma|f \right\rangle ~=~ \sum_{E} \left\langle E|\Gamma|f \right\rangle,
\]
where $E$ runs through all wide edges of $\Gamma$.

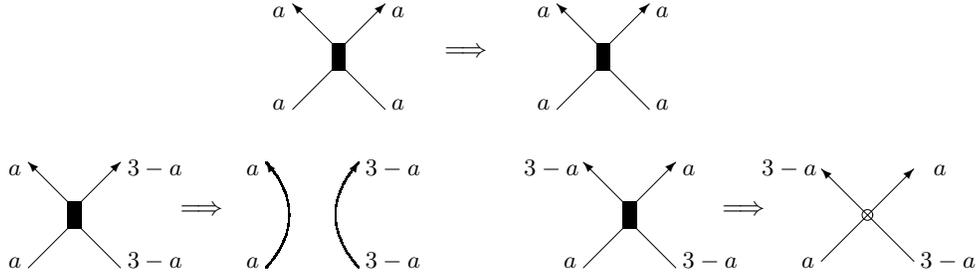
\begin{figure}[ht]

\setlength{\unitlength}{1pt}

\begin{picture}(360,100)(-180,0)


\put(-52.5,85){\vector(-1,1){15}}

\put(-47.5,85){\vector(1,1){15}}

\put(-52.5,75){\line(-1,-1){15}}

\put(-47.5,75){\line(1,-1){15}}

\put(-75,95){\small{$a$}}

\put(-75,60){\small{$a$}}

\put(-30,95){\small{$a$}}

\put(-30,60){\small{$a$}}

\put(-10,80){$\Longrightarrow$}


\put(52.5,85){\vector(1,1){15}}

\put(47.5,85){\vector(-1,1){15}}

\put(52.5,75){\line(1,-1){15}}

\put(47.5,75){\line(-1,-1){15}}

\put(70,95){\small{$a$}}

\put(70,60){\small{$a$}}

\put(25,95){\small{$a$}}

\put(25,60){\small{$a$}}


\put(-152.5,25){\vector(-1,1){15}}

\put(-147.5,25){\vector(1,1){15}}

\put(-152.5,15){\line(-1,-1){15}}

\put(-147.5,15){\line(1,-1){15}}

\put(-175,35){\small{$a$}}

\put(-175,0){\small{$a$}}

\put(-130,35){\small{$3-a$}}

\put(-130,0){\small{$3-a$}}

\put(-110,20){$\Longrightarrow$}


\qbezier(-77.5,0)(-60,20)(-77.5,40)

\put(-77.5,40){\vector(-2,3){0}}

\qbezier(-42.5,0)(-60,20)(-42.5,40)

\put(-42.5,40){\vector(2,3){0}}

\put(-85,35){\small{$a$}}

\put(-85,0){\small{$a$}}

\put(-40,35){\small{$3-a$}}

\put(-40,0){\small{$3-a$}}


\put(62.5,25){\vector(1,1){15}}

\put(57.5,25){\vector(-1,1){15}}

\put(62.5,15){\line(1,-1){15}}

\put(57.5,15){\line(-1,-1){15}}

\put(80,35){\small{$a$}}

\put(80,0){\small{$3-a$}}

\put(20,35){\small{$3-a$}}

\put(35,0){\small{$a$}}

\put(95,20){$\Longrightarrow$}


\put(167.5,2.5){\vector(-1,1){35}}

\put(132.5,2.5){\vector(1,1){35}}

\put(150,20){\circle{4}}

\put(175,35){\small{$a$}}

\put(170,0){\small{$3-a$}}

\put(110,35){\small{$3-a$}}

\put(125,0){\small{$a$}}

\linethickness{5pt}

\put(-50,75){\line(0,1){10}}

\put(50,75){\line(0,1){10}}

\put(-150,15){\line(0,1){10}}

\put(60,15){\line(0,1){10}}

\end{picture}

\caption{Modifying $\Gamma$}\label{label-change}

\end{figure}

Given a graph $\Gamma$ and a labeling $f\in\la(\Gamma)$, modify $\Gamma$ near each wide edge by the rules depicted in Figure \ref{label-change} to get a graph $\Gamma'$, where $a\in\{1,2\}$. All the regular edges of $\Gamma'$ labeled by $1$ (resp. 2) and wide edge adjacent to them form a closed graph $\Gamma_{f,1}$ (resp. $\Gamma_{f,2}$).

\begin{figure}[ht]

\setlength{\unitlength}{1pt}

\begin{picture}(360,40)(-180,60)


\put(-52.5,85){\vector(-1,1){15}}

\put(-47.5,85){\vector(1,1){15}}

\put(-52.5,75){\line(-1,-1){15}}

\put(-47.5,75){\line(1,-1){15}}

\put(-10,80){$\Longrightarrow$}


\qbezier(32.5,60)(50,80)(32.5,100)

\put(32.5,100){\vector(-2,3){0}}

\qbezier(67.5,60)(50,80)(67.5,100)

\put(67.5,100){\vector(2,3){0}}

\linethickness{5pt}

\put(-50,75){\line(0,1){10}}

\end{picture}

\caption{Removing a wide edge}\label{remove-wide-edge}

\end{figure}
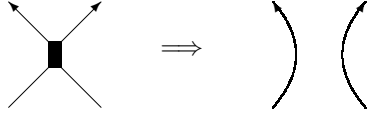

Also, for a closed graph $\Gamma$, removing all of its wide edges by the operation in Figure \ref{remove-wide-edge} gives a collection of disjoint oriented circles. Define the rotation number $r(\Gamma)$ of $\Gamma$ to be the total rotation number of this collection of circles.

For a graph $\Gamma$, a labeling $f\in\la(\Gamma)$ and $m,n\in \nat$, define
\[
\sigma_{m,n}(\Gamma,f) ~=~ \left\langle \Gamma|f \right\rangle + mr(\Gamma_{f,1})-nr(\Gamma_{f,2}).
\]

Wagner \cite{Wag} proved the following composition product formula.

\begin{theorem}[Theorem 1, \cite{Wag}]\label{wag-main}
For any closed graph $\Gamma$ and $i\in \zed$, $m,n\in \nat$,
\[
H_{m+n}^i(\Gamma) ~\cong~ \bigoplus_{f\in \la(\Gamma),~k+l+\sigma_{m,n}(\Gamma,f)=i} H_n^k(\Gamma_{f,1})\otimes_\C H_m^l(\Gamma_{f,2})\{\sigma_{m,n}(\Gamma,f)\}.
\]
\end{theorem}

Now we use the above theorem to prove following proposition, which implies Theorem \ref{res-braid}. (See \cite{Wu5} for the definition of resolved braids.)

\begin{proposition}[Proposition 3.8, \cite{Wu5}]\label{res-braid-prop}
Let $\Gamma$ be a resolved braid with $O$ strands and $e$ wide edges. For an integer $n\geq 1$, suppose that $\alpha$ is a non-zero homogeneous element of $H_n(\Gamma)$ (with respect to the quantum grading of $H_n(\Gamma)$). Then the quantum degree of $\alpha$ satisfies
\[
-(n-1)O-e ~\leq~ \deg \alpha ~\leq~ (n-1)O+e.
\]
\end{proposition}

\begin{proof}

We induct on the integer $n$. If $n=1$, as pointed out by Wagner \cite{Wag}, it is easy to see that
\[
H_1(\Gamma) ~\cong~
\left\{%
\begin{array}{ll}
    0, \text{ if } \Gamma \text{ has wide edges}, \\
    \C, \text{ if } \Gamma \text{ has no wide edges},  \\
\end{array}%
\right.
\]
where $\C$ has quantum grading $0$. So the proposition is trivially true for $n=1$.

Assume that the proposition is true for $H_n(\Gamma)$. Consider a non-zero homogeneous element $\alpha$ of $H_{n+1}(\Gamma)$. Let $m=1$ in the composition product formula. We get
\[
H_{n+1}^i(\Gamma) ~\cong~ \bigoplus_{f\in \la(\Gamma),~k+l+\sigma_{1,n}(\Gamma,f)=i} H_n^k(\Gamma_{f,1})\otimes_\C H_1^l(\Gamma_{f,2})\{\sigma_{1,n}(\Gamma,f)\}.
\]
So there is a labeling $f$ of $\Gamma$ and a non-zero homogeneous element $\beta$ of $H_{n}(\Gamma_{f,1})$ such that $\Gamma_{f,2}$ is a collection of disjoint circles, and
\[
\deg\alpha ~=~ \deg\beta+\sigma_{1,n}(\Gamma,f) ~=~ \deg\beta+ \left\langle \Gamma|f \right\rangle + r(\Gamma_{f,1})-nr(\Gamma_{f,2}).
\]
It is easy to see that
\[
r(\Gamma) ~=~ r(\Gamma_{f,1})+r(\Gamma_{f,2}),
\]
and
\[
-(e(\Gamma)-e(\Gamma_{f,1})) ~\leq~ \left\langle \Gamma|f \right\rangle ~\leq~ (e(\Gamma)-e(\Gamma_{f,1})).
\]
Also, since $\Gamma$, $\Gamma_{f,1}$ and $\Gamma_{f,2}$ are resolved braids in the same direction, one can see that $r(\Gamma)$, $r(\Gamma_{f,1})$ and $r(\Gamma_{f,2})$ have the same sign, and $O(\Gamma)=|r(\Gamma)|$, $O(\Gamma_{f,1})=|r(\Gamma_{f,1})|$ and $O(\Gamma_{f,2})=|r(\Gamma_{f,2})|$. By induction hypothesis, we have
\begin{eqnarray*}
\deg\alpha & = & \deg\beta+ \left\langle \Gamma|f \right\rangle + r(\Gamma_{f,1})-nr(\Gamma_{f,2}) \\
           & \leq & (n-1)O(\Gamma_{f,1}) + e(\Gamma_{f,1}) + (e(\Gamma)-e(\Gamma_{f,1})) + r(\Gamma_{f,1}) -n(r(\Gamma)-r(\Gamma_{f,1})) \\
           & = & nO(\Gamma)+e(\Gamma)+(n-1)(O(\Gamma_{f,1})+r(\Gamma_{f,1}))-n(O(\Gamma)+r(\Gamma)).
\end{eqnarray*}
But it is clear that $O(\Gamma)+r(\Gamma) \geq O(\Gamma_{f,1})+r(\Gamma_{f,1}) \geq 0$. Thus,
\[
\deg\alpha ~\leq~ nO(\Gamma)+e(\Gamma).
\]
Similarly,
\begin{eqnarray*}
\deg\alpha & = & \deg\beta+ \left\langle \Gamma|f \right\rangle + r(\Gamma_{f,1})-nr(\Gamma_{f,2}) \\
           & \geq & -(n-1)O(\Gamma_{f,1}) - e(\Gamma_{f,1}) - (e(\Gamma)-e(\Gamma_{f,1})) + r(\Gamma_{f,1}) -n(r(\Gamma)-r(\Gamma_{f,1})) \\
           & = & -nO(\Gamma)-e(\Gamma)-(n-1)(O(\Gamma_{f,1})-r(\Gamma_{f,1}))+n(O(\Gamma)-r(\Gamma)) \\
           & \geq & -nO(\Gamma)-e(\Gamma),
\end{eqnarray*}
since $O(\Gamma)-r(\Gamma) \geq O(\Gamma_{f,1})-r(\Gamma_{f,1}) \geq 0$. Thus, the proposition is true for $H_{n+1}(\Gamma)$. This completes the proof.
\end{proof}

\begin{proof}[Proof of Theorem \ref{res-braid}]
The proof is now very easy, and is identical to that in \cite{Wu5}. We sketch it here for the convenience of the reader.

$H_n(B)$ is the cohomology of a homogeneous differential map on the graded $\C$-linear space
\[
\bigoplus_{\Gamma} H_n(\Gamma)\{(n-1)w(B)+e_+(\Gamma)-e_-(\Gamma)\},
\]
where $\Gamma$ runs through all resolutions of $B$, and $e_+(\Gamma)$ (resp. $e_-(\Gamma)$) is the number of wide edges in $\Gamma$ from resolving positive (resp. negative) crossings of $B$. Then Theorem \ref{res-braid} from Proposition \ref{res-braid-prop} and the fact that $e_\pm(\Gamma)\leq c_\pm(B)$.
\end{proof}

\begin{remark}
Wagner's composition product for the Khovanov-Rozansky cohomology is a generalization of Jaeger's composition product for the HOMFLY polynomial \cite{Jaeger}. Jaeger \cite{Jaeger} also gave an alternative proof of inequality \eqref{HOMFLYBi}. Our proof here can be considered a generalization of Jaeger's proof.
\end{remark}

\section{Proofs of Theorems \ref{bennequin} and \ref{sharper-bennequin}}

The main technical tool in this section is link cobodism. In particular, we will use the next proposition repeatedly.

\begin{proposition}\label{isomorphism}
If $S$ is a connected cobodism from knot $K_1$ to knot $K_2$, then $S$ induces a $\C$-linear isomorphism $\Psi_S:H_p(K_1)\rightarrow H_p(K_2)$ of quantum degree $\leq -(n-1)\chi(S)$.
\end{proposition}
\begin{proof}
First find a movie presentation of $S$. This movie presentation induces a $\C$-linear homomorphism $\Psi_S:H_p(K_1)\rightarrow H_p(K_2)$ of quantum degree $\leq -(n-1)\chi(S)$. From Proposition 5.1 of \cite{Wu7}, one can see that $\Psi_S$ is surjective. But, by Theorem 2 of \cite{Gornik}, $\dim_\C H_p(K_1) = \dim_\C H_p(K_1)=n$. So $\Psi_S$ is in fact an isomorphism.
\end{proof}

\begin{remark}
Note that $\Psi_S$ is defined using a movie present. It is not clear whether $\Psi_S$ is independent (up to scaling) of the choice of the movie presentation. But this does not affect our proofs.
\end{remark}

\begin{proof}[Proof of Theorem \ref{bennequin}]
From the Corollary in Section 4 of \cite{Ya}, we only need to prove this theorem for closed braid diagrams of knots.

Assume that $B$ is a closed braid diagram for a knot with only positive crossings. Then, by Theorems \ref{slice-genus} and \ref{res-braid} of the present paper and Theorem 1.2 of \cite{Wu7}, we have
\[
(n-1)(w(B)-O(B)) ~\leq~ \gnm(B) ~\leq~ \gpm(B) ~\leq~ -(n-1)\chi_s(B) ~\leq~ (n-1)(w(B)-O(B)),
\]
where the last inequality is true because that $O(B)-w(B)$ is the Euler characteristic of the Seifert surface of $B$ obtained from the Seifert algorithm. Thus
\[
\gpm(B)=(n-1)(w(B)-O(B)).
\]

Next we prove by inducting on the number of negative crossings that, for any closed braid diagram $B$ of a knot,
\begin{equation}\label{bennequin-min}
(n-1)(w(B)-O(B)) ~\leq~ \gpm(B).
\end{equation}
If $B$ has no negative crossings, then Inequality
\eqref{bennequin-min} is true from the above equation. Assume that
$B$ is a closed braid diagram of a knot with negative crossings, and
\eqref{bennequin-min} is true for any closed braid diagrams of knots
with less negative crossings than $B$. Let $c$ be a negative
crossing of $B$, and $B^+$ the closed braid obtained from $B$ by
changing $c$ into a positive crossing. Consider the link cobodism
$S$ from $B$ to $B^+$ with the movie presentation in Figure
\ref{movie1}. Then $\chi(S)=-2$. Let $\alpha$ be an non-zero element
of $H_p(B)$ with quantum degree $\gpm(B)$. Then, since $\Psi_S$ is
an isomorphism, $\Psi_S(\alpha)\neq0$, and, therefore, has quantum
degree at least $\gpm(B^+)$. By induction hypothesis and Proposition
\ref{isomorphism},
\[
(n-1)(w(B^+)-O(B^+)) ~\leq~ \gpm(B^+) ~\leq~ \deg (\Psi_S(\alpha)) ~\leq~ \gpm(B) + 2(n-1).
\]
Note that $w(B^+)=w(B)+2$ and $O(B^+)=O(B)$. This implies that \eqref{bennequin-min} is true for $B$. Thus, we have proved Inequality \eqref{bennequin-min}.

\begin{figure}[ht]

\setlength{\unitlength}{1pt}

\begin{picture}(360,50)(-180,-30)


\put(130,-20){\vector(1,1){40}}

\put(170,-20){\line(-1,1){15}}

\put(145,5){\vector(-1,1){15}}

\put(145,-30){$B^+$}


\qbezier(30,-20)(50,0)(30,20)

\put(30,20){\vector(-1,1){0}}

\put(70,20){\vector(2,3){0}}

\put(55,-5){\line(1,-1){15}}

\qbezier(70,20)(45,-20)(45,0)

\qbezier(45,0)(45,5)(51,-1)

\put(90,0){$\Longrightarrow$}


\qbezier(-70,-20)(-50,0)(-70,20)

\put(-70,20){\vector(-1,1){0}}

\put(-45,5){\vector(1,1){15}}

\qbezier(-30,-20)(-55,20)(-55,0)

\qbezier(-55,0)(-55,-5)(-49,1)

\put(-10,0){$\Longrightarrow$}


\put(-130,-20){\vector(-1,1){40}}

\put(-170,-20){\line(1,1){15}}

\put(-145,5){\vector(1,1){15}}

\put(-155,-30){$B$}

\put(-110,0){$\Longrightarrow$}

\end{picture}

\caption{the cobodism $S$}\label{movie1}

\end{figure}
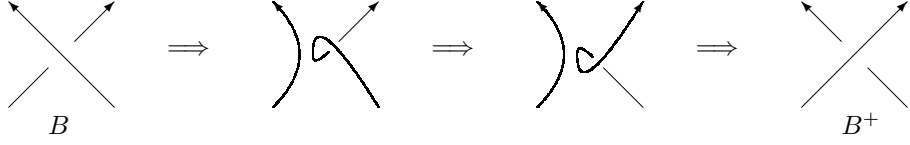

Similarly, one can prove that, for any closed braid diagram $B$ of a knot,
\begin{equation}\label{bennequin-max}
\gpM(B) ~\leq~ (n-1)(w(B)+O(B)).
\end{equation}
To do that, first show that $\gpM(B)=(n-1)(w(B)+O(B))$ for a closed braid diagram $B$ of a knot with only negative crossings. Then prove \eqref{bennequin-max} by inducting on the number of positive crossings. The details are left to the reader.
\end{proof}

\begin{proof}[Proof of Theorem \ref{sharper-bennequin}] (a) Assume $D_K$ has only positive crossings. By Theorems \ref{slice-genus} and \ref{bennequin},
\[
(n-1)(w-O) ~\leq~ \gpm(K) ~\leq~ -(n-1)\chi_s(K) ~\leq~ (n-1)(w-O),
\]
where, again, the last inequality is true since $O-w$ is the Euler characteristic of the Seifert surface of $K$ obtained by applying the Seifert Algorithm to $D_K$. Thus, $\gpm(K)=-(n-1)\chi_s(K)=(n-1)(w-O)$. Moreover, by Theorem \ref{slice-genus},
\[
\gpM(K)\leq(n-1)(2-\chi_s(K))= (n-1)(2+w-O).
\]
So we have proved part (a). The proof of part (b) is similar, and is left to the reader.

(c) We apply an argument by Rudolph \cite{Ru1} to prove part (c). This method was also used by Kawamura \cite{Kawa} to establish \eqref{KBi}. We only show that
\begin{equation}
\gpM(K) ~\leq~ (n-1)(w+O_\leq-O_>),
\end{equation}
and leave the proof of the other half of the inequality to the reader.

Let $c_+$ (resp. $c_-$) be the number of positive (resp. negative) crossings in $D_K$. Apply the Seifert algorithm to every positive crossing of $D_K$, which modify $D_K$ into a disjoint union of a negative link diagram $L^-$ and a collection of $O_>$ disjoint circles. There is a link cobodism $S_1$ from $L^-$ to $K$ with Euler characteristic $O_>-c_+$ consisting of $O_>$ $0$-handles (circle creations), $c_+$ Reidemeister I moves creating $c_+$ positive kinks, and $c_+$ $1$-handles (saddle moves).

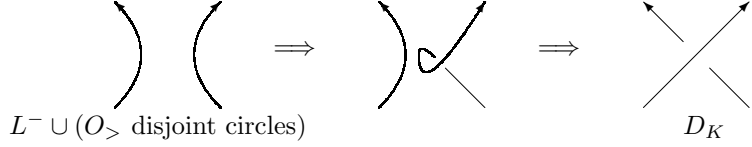
\begin{figure}[ht]

\setlength{\unitlength}{1pt}

\begin{picture}(360,50)(-180,-30)


\put(80,-20){\vector(1,1){40}}

\put(120,-20){\line(-1,1){15}}

\put(95,5){\vector(-1,1){15}}

\put(95,-30){$D_K$}


\qbezier(-20,-20)(0,0)(-20,20)

\put(-20,20){\vector(-1,1){0}}

\put(20,20){\vector(2,3){0}}

\put(5,-5){\line(1,-1){15}}

\qbezier(20,20)(-5,-20)(-5,0)

\qbezier(-5,0)(-5,5)(1,-1)

\put(40,0){$\Longrightarrow$}


\qbezier(-120,-20)(-100,0)(-120,20)

\put(-120,20){\vector(-1,1){0}}

\put(-80,20){\vector(1,1){0}}

\qbezier(-80,-20)(-100,0)(-80,20)

\put(-60,0){$\Longrightarrow$}

\put(-160,-30){$L^-\cup(O_>$ disjoint circles$)$}

\end{picture}

\caption{A $1$-handle of the cobodism $S_1$}\label{movie2}

\end{figure}

Let $k$ be the number of components of the link represented by $L^-$. We joint all these components by $k-1$ movies of type (i) or (ii) depicted in Figure \ref{movie3} by following steps:
\begin{enumerate}[(I)]
    \item Let $D_1, \cdots,D_l$ be the components of $L^-$ as a subset of the plane. Use $l-1$ movies of type (i) or (ii) to connect these components. This results in a link diagram $L'$ with only negative crossings which is a connected subset of the plane.
    \item Use movies of type (ii) to connect the components of the link represented by $L'$ resulting in a diagram $K^-$ of a knot with only negative crossings.
\end{enumerate}
It is clear that a movie of type (ii) reduces the writhe by $1$, and does not change the number of Seifert circles of the diagram. Also, movies of type (i) do not change the writhe. Moreover, since we only applied movies of type (i) to connect between different components of $L^-$ as a subset of the plane, every type (i) movie we used reduces the number of Seifert circles by $1$. Let $m_1$ and $m_2$ be the numbers of type (i) and type (ii) movies used in the above construction. Then
\begin{eqnarray*}
k-1 & = & m_1+m_2, \\
w(K^-)& = & w(L^-)-m_2 ~=~ -c_--m_2, \\
O(K^-) & = & O(L^-)-m_1 ~=~ O_\leq -m_1.
\end{eqnarray*}

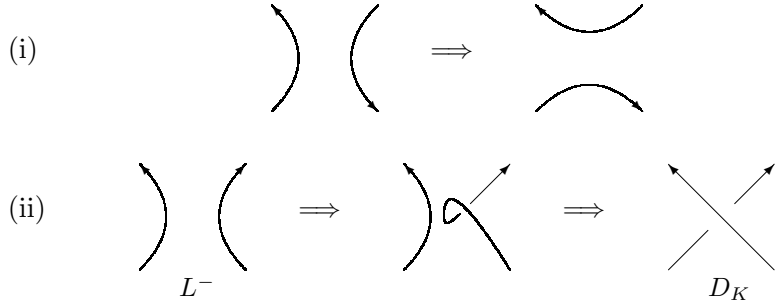
\begin{figure}[ht]

\setlength{\unitlength}{1pt}

\begin{picture}(360,110)(-180,-30)


\put(-170,60){\large{(i)}}

\qbezier(-70,40)(-50,60)(-70,80)

\put(-70,80){\vector(-1,1){0}}

\qbezier(-30,40)(-50,60)(-30,80)

\put(-30,40){\vector(1,-1){0}}

\put(-10,60){$\Longrightarrow$}


\qbezier(70,40)(50,60)(30,40)

\put(30,80){\vector(-1,1){0}}

\qbezier(70,80)(50,60)(30,80)

\put(70,40){\vector(1,-1){0}}


\put(-170,0){\large{(ii)}}

\qbezier(-120,-20)(-100,0)(-120,20)

\put(-120,20){\vector(-1,1){0}}

\put(-80,20){\vector(1,1){0}}

\qbezier(-80,-20)(-100,0)(-80,20)

\put(-60,0){$\Longrightarrow$}

\put(-105,-30){$L^-$}


\qbezier(-20,-20)(0,0)(-20,20)

\put(-20,20){\vector(-1,1){0}}

\put(5,5){\vector(1,1){15}}

\qbezier(20,-20)(-5,20)(-5,0)

\qbezier(-5,0)(-5,-5)(1,1)

\put(40,0){$\Longrightarrow$}


\put(120,-20){\vector(-1,1){40}}

\put(80,-20){\line(1,1){15}}

\put(105,5){\vector(1,1){15}}

\put(95,-30){$D_K$}

\end{picture}

\caption{Movies (i) and (ii)}\label{movie3}

\end{figure}

Reversing these movies, we get a cobodism $S_2$ from $K^-$ to $L^-$ with Euler characteristic $1-k$. Consider the cobodism $S=S_2\cup_{L^-}S_1$ from $K^-$ to $K$. It has Euler characteristic $O_>-c_+-(k-1)$, and, by Proposition \ref{isomorphism}, induces an isomorphism
\[
\Psi_S:H_p(K^-)\xrightarrow{\cong} H_p(K)
\]
of quantum degree $\leq (n-1)(c_+-O_>+k-1)$. By Part (b) of the theorem, $\gpM(K^-)=(n-1)(w(K^-)+O(K^-)$. Let $\alpha$ be a non-zero element of $H_p(K)$ of quantum degree $\gpM(K)$. Then there exists a non-zero element $\beta$ of $H_p(K^-)$ such that $\Psi_S(\beta)=\alpha$. So
\begin{eqnarray*}
\gpM(K) & = & \deg \alpha \\
        & \leq & \deg\beta + (n-1)(c_+-O_>+k-1) \\
        & \leq & \gpM(K^-) + (n-1)(c_+-O_>+k-1) \\
        & = & (n-1)(w(K^-)+O(K^-)+c_+-O_>+k-1) \\
        & = & (n-1)(w+O_\leq-O_>).
\end{eqnarray*}
\end{proof}

\end{document}